\documentclass[12pt,a4paper]{article}
\usepackage{amsmath,amssymb,amsthm,amsfonts,amscd,euscript,verbatim}
\usepackage{t1enc}
 \theoremstyle{plain}
\newtheorem{theo}{Theorem}[section]

\newtheorem{pr}[theo]{Proposition}

\theoremstyle{remark}
\newtheorem{rema}[theo]{Remark}

\theoremstyle{definition}

\newcommand\smc{{SmCor}}

\newcommand\dmge{DM^{eff}_{gm}{}}
\newcommand\dmgeq{DM^{eff}_{gm}\q{}}
\newcommand\dmgm{DM_{gm}}

\newcommand\mg{M_{gm}}
\newcommand\obj{Obj}

\newcommand\q{{\mathbb{Q}}}
\newcommand\af{\mathbb{A}}

\newcommand\chow{Chow}
\newcommand\chowe{Chow^{eff}}

 \DeclareMathOperator\ke{\operatorname{Ker}}
 \DeclareMathOperator\cok{\operatorname{Coker}}

\begin{document}
 \title{Explicit  generators for (conjectural) mixed motives (in  Voevodsky's $\dmge$). The  Kunneth decomposition of pure (numerical) motives }
 \author{M.V. Bondarko
 }
 \maketitle


 \section*{Introduction}

The current version of this text is just an announcement of
results. It includes the main ideas that are quite sufficient to
prove the  results announced; yet complete proofs will be written
down later.

We recall that constructing the category of mixed motives $MM$ is
probably the most challenging problem of the modern algebraic
geometry. Note that $MM$ should be a full subcategory  of
Voevodksy's $\dmge$ or $\dmge\q$. Yet Voevodsky's motives  that
are currently assumed to be mixed do not generate $\dmge$; they do
not generate even motives of dimension $\le 2$ (to the knowledge
of the author).

In this note we describe very explicitly a rich family of mixed
motives that surely generates $\dmgeq$ (as a triangulated
category). They "should be" mixed since they have only one
non-zero Betti cohomology group. Our method also allows to define
a family of direct summands of the numerical motif of any smooth
projective variety $P$. Modulo certain standard conjectures, this
construction yields the Kunneth decomposition of the diagonal of
$P$.

From the "motivic" point of view the main ideas of this paper are:
repetitive hyperplane sections of a smooth affine variety give a
resolution of its Voevodksy's motif by mixed motives; applying
this fact to a Jouanalou's "replacement" of a smooth projective
variety $P$ one obtains  the Kunneth decomposition of the diagonal
of the numerical motif of $P$.

It seems that our ideas are related to the motivic ideas of M.
Nori. Yet the results described in \S2 seem to be completely new.

In this version we will assume that our base field $k$ is fixed
and  canonically embedded in the field of complex numbers. We will
denote by $H_B$ the (rational) Betty cohomology of a variety or a
motif.  We will only consider  motives and cohomology with
rational coefficients.

We will use some definitions and notation of \cite{1}. $\dmge\q$
will denote the idempotent completion of $\dmge$ with morphisms
tensored by $\q$. Similarly, $\chow$ ($\chowe$) will denote the
category of rational (effective) Chow motives; $Mot_{hom}$ will be
the category of rational effective homological motives;
$Mot_{num}$ will be the category of effective rational numerical
motives.

We will also  mention the (conjectural) abelian category of
(effective) mixed motives $MM\subset \dmgeq$; so the (conjectural)
mixed motivic cohomology functors will be covariant.
 $H^i_M(X)\in\obj MM$ will denote   the
$i$-th cohomology of  $X\in\obj\dmgeq$ with respect to the
(conjectural!) mixed motivic structure.  Often $X$ will be the
rational Voevodsky's motif of a (smooth) variety $V$ (denoted by
$\mg(V)$).


\section{A method for constructing mixed motives (in $\dmge$)}

The starting point of this work were certain results of \cite{a1}
and \cite{a2}.

Let $U$ be a smooth Zarisky open subvariety of a smooth projective
variety $Y$ of dimension $m$ and let $H$ be its hyperplane
section; let $V=U\cap H$ be non-empty.

We recall Theorem 6.1.1 of \cite{a1}. It states that the induced
map $H^i(V)\to H^i(U)$ is bijective for $i<m-1$ and is injective
for $i=m-1$. To the author's knowledge this means that $(U,V)$ is
a {\it good pair} in the sense of M. Nori.

\begin{pr}\label{fir}
Let $U$ be affine. We define $M\in\obj \dmgeq$ as a cone of
$\mg(V)\to \mg(U)$ i.e. $\mg(V)\to \mg(U)\to M$ is a distinguished
triangle. Then the only non-zero Betty cohomology group of $M$ is
$H_B^m$.
\end{pr}
\begin{proof}
For any $i$ we have a long exact sequence $$\dots\to  H^{i-1}_B(U)
\to H^{i-1}_B(V) \to H^i_B(M) \to H^i_B(U) \to H^i_B(V)\to
\dots.$$

It remains to recall that $H^j_B(X)=0$ for $j>m$ and apply Theorem
6.1.1 of \cite{a1}.

\end{proof}

\begin{rema}
1.  It follows that $M[-m]$ a nice candidate for an (abelian)
mixed motif. We will call  motives of this type "our" motives. For
example, the becomes mixed in Hanamura's construction (see \S3 of
\cite{ha3}; note that Hanamura's construction is purely
conjectural!). Recall that Hanamuras's motivic category is
anti-equivalent to (Voevodskys's) $\dmgm\q$; see the proof in \S4
of \cite{mymot}. If $M[-m]$ in $MM$  should be an extension of
$H_M^m(X)$ by $\ke:H_M^{m-1}(Y)\to H_M^{m-1}(X)$. Note that the
corresponding fact holds in the category of (rational) mixed Hodge
structures i.e.  $H^i(M)_{Hodge}$ is an extension of
 $\cok:H_{Hodge}^{m-1}(X)\to
H_{Hodge}^{m-1}(Y)$ by $H_{Hodge}^m(X)$.

2. One could easily construct a large variety of "$1$-extensions"
of our motives. This leads to the following construction of mixed
motives that is  (at the moment) the most general among those
known to the author.

We take a complex $X_i\in K^b(\smc\otimes\q)$ with $\dim (X_i)=i$;
$X_i$ are affine. Let $Y_i\subset X_i$ be a compatible system of
hyperplane sections (as in Proposition \ref{fir}). Then the cone
of $(Y_i)\to (X_i)$ (considered as a motif coming from
$K^b(\smc\otimes\q)$) cannot have any cohomology out of dimension
$0$. This statement follows from Proposition \ref{fir} by easy
induction. In particular, in this way one could describe tensor
products of "our" motives in $\dmgeq$. (Recall that the tensor
product of motives is compatible with the tensor product of
varieties.) 

It also seems that the "twisted dual" of $M$ i.e a cone $M'$ of
the natural morphism $\mg^c(U)\to \mg^c(V)(1)[2] $ can be obtained
from $K^b(\smc\otimes\q)$ using this construction (see Theorem
4.3.7 of \cite{1}). One could also easily check that $M'$ has only
one non-zero cohomology group using the Poincare duality and
Proposition \ref{fir}.

3. It is easily seen that "our" motives  generate the whole
$\dmgeq$ as a triangulated category (after idempotent completion).
Indeed, a repetitive application of Proposition \ref{fir} (to $U$,
then to $V$, etc.) immediately yields that the triangulated
category $D$ generated by "our" motives contains all motives of
smooth affine varieties. This method should correspond to the
resolution of $\mg(U)$ by mixed motives.

Next, the Mayer-Viertoris triangle (see \S2 of \cite{1}) yields
that $D$ contains motives of all smooth varieties. Lastly, recall
that $\dmgeq$ is the idempotent completion of a certain
localization of $K^b(\smc\otimes\q)$.

 Moreover, it seems that any object of $\dmgeq$ that
comes from $K^b(\smc\otimes\q)$ could be presented as a certain
"complex" of "our" motives  (i.e. it becomes equivalent in
$\dmgeq$ to such a complex constructed inside of
$K^b(\smc\otimes\q)$).

4. One could try to define the mixed motivic $t$-structure using
"our" motives. Then is seems that the main problem is to prove
that there are no morphisms of negative degrees between "our"
motives (in $\dmgeq$). Unfortunately,  one probably cannot
overcome this difficulty with assuming certain vanishing
conjectures (as it was done in \cite{ha3}). Still there is a hope
to deduce everything from the (conjectural) conservativity of the
Betti realization of motives. Note that (by the conservativity of
the weight complex functor, see Proposition 6.1.3 of \cite{mymot},
and by the existence of the weight spectral sequence for
realizations of motives, see (14) of \S7.3 ibid.) it suffices to
check this conservativity on $K^b(\chowe)$.

\end{rema}


It seems that Proposition \ref{fir} is should be connected with
the ioga of relative motives (mainly of motives over $\q(i)$).

\section{A candidate for the Kunneth decomposition of pure
motives}

Unfortunately, it is not clear how to describe all $H^i_M(U)$ for
all $i$ using our method. Yet this seems to be easy if we restrict
ourselves to numerical motives (of smooth projective varieties).

More precisely,  we apply Theorem 2.1 of \cite{a2} for smooth
projective $X=Y$ (in the notation of loc.cit.). Essentially, this
is (more or less) equivalent to the composition of a repetitive
application of Proposition \ref{fir} with Jouanolou's trick and
Theorem 6.1.1 of \cite{a1}. In order to clarify the connection of
the reasoning below with Proposition \ref{fir} we recall that for
any smooth quasi-projective $X$ there exists an affine line bundle
(of some dimension $l\ge 0$) $X'/X$ such that $X'$ is affine over
$k$ (Jouanolou's trick). We will call $X'$ a {\it Jouanolou's
replacement} for $X$. Since $X'$ is Zarisky locally isomorphic to
$X\times \af^l$, we obtain that $\mg(X')\cong \mg(X)$ in
$\dmge\q$; cf. Proposition 3.5.1 of \cite{1}.


Now, Theorem 2.1 of \cite{a2} applied to smooth projective $X=Y$
(in the notation of loc.cit.) states exactly that the canonical
filtration for the cohomology of $X$ can be described by kernels
of $H_B(X)\to H_B(X_p)$ for some morphism of smooth varieties
$X_p\to X$ (here $X_p$ are obtained by taking smooth hyperplane
sections of any Jouanolou's replacement $X'$ of $X$). Now one can
look at the numerical motif of $X$. $X_p$ is (probably) not
projective; yet it corresponds to a complex of Chow motives (see
\S6 of \cite{mymot} for the consideration of {\it weight complex}
of motives). Now consider
the cokernel 
 of the maps induced by $X_p\to X$ on the
(covariant) numerical motif of $X$. More precisely, by part 1 of
Theorem 6.2.1 of \cite{mymot} the weight complex $t(\mg(X_p))$ is
a complex of (effective) Chow motives concentrated in non-negative
degrees. We consider the corresponding complex
$T_p=t_{num\q}(\mg(X_p))=T_0\to T_1\to \dots \in K^b(Mot_{num})$.
Note that $X_p$ could be taken to be a complement of one smooth projective varieties
by another one; hence we may assume that $T_i=0$ for $i>1$. 
 Since the category of rational (effective)
numerical motives is abelian semi-simple (see Theorem 1 of
\cite{ja}), the cokernel of $H^0(T_p)\to Mot_{num\q}(X)$ gives a
direct summand of the motif of $X$. It seems to be a very
reasonable candidate for $\sum_{p\le i\le 2n}H_M^i(X)[-i]$ (as a
numerical motif with the grading compatible with the theory of
mixed motives). Note that $X_p$ are embedded into each other, so
one can define each $H_M^i(X)$ separately (unconditionally). We
obtain an explicit construction that (conjecturally) yields the
Kunneth decomposition of the diagonal of the numerical motif of
$X$. Indeed, the construction obviously gives the the Kunneth
decomposition if the numerical equivalence is equivalent to the
homological equivalence. Note also that (using easy finite
dimension arguments similar to those of \cite{ja}) any direct
summand of the numerical motif of $X$ could be lifted to a direct
summand of its homological motif (for example, one could consider
the "minimal possible lift"). Yet it doesn't seem easy to prove
that the homological motives obtained give the Kunneth
decomposition on this level. Still even without this fact the
decomposition constructed could help to describe the (conjectural)
Tannakian category of pure motives.

It doesn't seem to be very difficult to prove (certain)
functoriality of this construction; in particular, to verify that
the construction is canonical. Possibly, to simplify the proof
one should modify the construction a little. Probably this leads
to the consideration of (a certain) coniveau filtration of the
motif of $X$. Note that to this end we should consider some
Jouanolou's replacement  $X'$ of $X$.  Possibly it also makes
sense to consider a certain "limit" with respect to all
Jouanolou's replacements for $X$; note the fiber product over $X$
of two Jouanolou's replacements for $X$ is also a  replacement for
it.

Lastly, we note that the construction could also be applied on the
level of Chow motives. It does not necessarily yield a Chow motif
since $\chowe$ is not abelian; yet it yields a factor of the
functor represented by $Mot_{Chow}(X)$ in the abelian category
$\chowe_*$ of additive contravariant functors $\chowe\to AbGr$ (we
have a natural full embedding $\chowe\to \chowe_*$). It doesn't
seem difficult to make this construction functorial (using the
limit methods described above). Note that then our construction
would automatically give canonical functors $F^i_{Hom}:\chowe\to
Mot_{hom*}$ and $F^i_{num}:\chowe\to Mot_{num}$ (since $Mot_{num}$
is abelian!). This approach could also yield (by applying Theorem
2.1 of \cite{a2}) that these functors could be constructed using
the 'basic' construction (without passing to any limits); the
corresponding "Chow fact" seems to be more difficult.

A curious observation:  for any $X_p$ all $T_i$ (could be chosen
to be) divisible by $\q(i)$. Hence for the birational motives
theory (where $\q(1)$ is killed) it suffices to consider only the
cokernel of $T_0\to Mot_{num}(X)$; see \cite{kahnbi}.

\section{Concluding remarks}

Certainly, these arguments cannot prove all standard conjectures.
Yet they are very explicit; this gives a hope to prove some parts
of standard conjectures (or other interesting statements)
unconditionally. Also, there should be a considerable impact on
interrelations between conjectures. The author's knowledge of this
field (of conjectures) is quite poor; so he will be deeply
grateful for any ideas in this direction.


\begin{thebibliography}{1}

\bibitem{a1}
 Intermediate Jacobians and Hodge Structures of Moduli
Spaces, Arapura D., Sastry P. // Proc. Indian Acad.  Sci. (Math
Sci.), vol. 110, No. 1, 2000, 1--26;
 see
http://www.ias.ac.in/mathsci/vol110/feb2000/Pm1709.pdf


\bibitem{a2}
The Leray spectral sequence is motivic, Arapura D.// Inv. Math.,
vol. 160 (2005), 567--589;
 see also
http://arxiv.org/abs/math.AG/0301140.


\bibitem{mymot} Bondarko M.V.,
{'Weight enhancement' and truncations for Voeovodsky motives;
filtrations and motivic descent spectral sequence  for
differential graded
 realizations}, electronic,
 http://arxiv.org/abs/math.AG/0601713.

\bibitem{ha3} Hanamura M. Mixed motives and algebraic cycles, III//
Math. Res. Letters 6, 1999,  61--82.



\bibitem{ja} Jannsen U., Motives, numerical equivalence, and semi-simplicity //Inv. Math.,
vol. 107 (1992), 447--452.


\bibitem{kahnbi} Kahn B., Sujatha R., Birational motives, I,
preprint, available on the K-theory preprint archive, no 0596,
http://www.math.uiuc.edu/K-theory/0596/

  \bibitem{1} Voevodsky V. Triangulated category of motives, in:
 Voevodsky V.,  Suslin A., and  Friedlander E.,
Cycles, transfers and motivic homology theories, Annals of
 Mathematical studies, vol. 143, Princeton University Press,
 2000,  188--238,  see also http://www.math.uiuc.edu/K-theory/0074/


\end{thebibliography}
\end{document}